\documentclass[fleqn,10pt]{article}
\usepackage{listings}
\usepackage{xcolor}
\definecolor{codegreen}{rgb}{0,0.6,0}
\definecolor{codegray}{rgb}{0.5,0.5,0.5}
\definecolor{codepurple}{rgb}{0.58,0,0.82}
\definecolor{backcolour}{rgb}{0.95,0.95,0.92}
\usepackage{caption}
\usepackage{subcaption}

\usepackage{mathtools,cases,csquotes}
\usepackage{cancel}
\usepackage{amsbsy,csquotes}
\usepackage{booktabs}
\usepackage{subcaption}
\usepackage{amsmath, amssymb}
\usepackage{listings}
\usepackage{empheq}
\usepackage{afterpage}
\usepackage{xcolor}
\usepackage{cleveref}
\usepackage{url}
\usepackage{booktabs}
\usepackage{stmaryrd}
\usepackage{color}
\usepackage{graphicx}
\usepackage{overpic}
\usepackage{subcaption}

\lstdefinestyle{mystyle}{
    backgroundcolor=\color{backcolour},   
    commentstyle=\color{codegreen},
    keywordstyle=\color{magenta},
    numberstyle=\tiny\color{codegray},
    stringstyle=\color{codepurple},
    basicstyle=\ttfamily\footnotesize,
    breakatwhitespace=false,         
    breaklines=true,                 
    captionpos=b,                    
    keepspaces=true,                 
    numbers=left,                    
    numbersep=5pt,                  
    showspaces=false,                
    showstringspaces=false,
    showtabs=false,                  
    tabsize=2
}
\usepackage[a4paper, total={5.5in, 8in}]{geometry}

\title{\textsc{graphnics}: Combining \textsc{FEniCS} and \textsc{NetworkX} to simulate flow in complex networks}

\author{Ingeborg G. Gjerde}

\begin{document}

\flushbottom
\maketitle
\thispagestyle{empty}

\section*{Summary}
Network models facilitate inexpensive simulations, but require careful handling of bifurcation conditions. We here present the \textsc{graphnics} library \cite{graphnics}, which combines \textsc{FEniCS} \cite{logg2012automated} with \textsc{NetworkX} \cite{hagberg2008exploring} to facilitate network simulations using the finite element method. \textsc{graphnics} features
\begin{itemize}
    \item A \texttt{FenicsGraph} class built on top of the \textsc{NetworkX} \texttt{DiGraph} class, that constructs a global mesh for a network and provides \textsc{fenics} mesh functions describing how they relate to the graph structure. 
    \item Example models showing how the \texttt{FenicsGraph} class can be used to assemble and solve different network flow models.
\item Demos showing e.g. how the simulations in \cite{daversin2022geometrically} can be extended to complex biological networks. Interestingly, the results show that vasomotion modelled as a travelling sinusoidal wave is capable of driving net perivascular fluid flow through an arterial tree, as was proposed in \cite{van2020vasomotion} based on experimental data.
\end{itemize}
The example models are implemented using \textsc{fenics}$\_$\textsc{ii} \cite{kuchta2021assembly}. Only minor adaptions are needed to adapt the code to use e.g. the mixed-dimensional branch of \textsc{FEniCS} \cite{daversincatty2019abstractions}, as the assembly uses the common block structure of the problem. 

\section*{Statement of need}
\textsc{FEniCS} \cite{logg2012automated} provides high-level functionality for specifying variational forms. This allows the user to focus on the model they are solving rather than implementational details. Network models typically impose conservation of mass at bifurcation points, i.e. that there should be no jump in the cross-section flux over junctions. \textsc{FEniCS} implicitly assumes that each vertex is connected to two cells. At bifurcation vertices in a network, however, the vertex is connected to three or more cells. Thus one cannot use the jump operator currently offered in \textsc{FEniCS}. Moreover, the manual assembly of the jump terms is highly prone to errors once the network becomes non-trivial.

In this software we extend the \texttt{DiGraph} class offered by \textsc{NetworkX} so that it (i) creates a global mesh of the network and (ii) creates data structures describing how this mesh is connected to the graph structure of the problem. In addition to this we provide convenience functions that then make it straightforward to assemble and solve network flow models using the finite element method.

This is showcased in a series of demos; in particular, we show how to extend the network flow simulations in \cite{daversin2022geometrically} so that they run on complex vascular domains. The simulations consider perivascular fluid flow driven by arterial wall pulsations. Experimental results indicate that artertial wall motion drives bulk perivascular fluid flow \cite{mestre2018flow}. Simulation studies have so far been unable to recreate this \cite{daversin2020mechanisms,daversin2022geometrically}, perhaps due to a lack of network complexity. Using the \textsc{NetworkGen} library \cite{networkgen} to generate vascular trees, we find that vasomotion can drive bulk fluid flow in arterial trees with several branching generations.

\section*{Mathematics}

Let $\mathsf{G}=(\mathsf{V}, \mathsf{E})$ denote a graph $\mathsf{G}$ with $n$ edges $\mathsf{E}$ and $m$ vertices $\mathsf{V}$. We let $\Lambda_i$ be the geometrical domain associated with the edge $e_i \in \mathsf{E}$. 

We want to solve flow models on networks, one example being the hydraulic network model
\begin{align}
    \mathcal{R} \hat{q}_i - \partial_s \hat{p}_i &= 0 \text{ on } \Lambda_i, && \text{(constitutive equation on edge)} \label{eq:1}\\ 
    \partial_s \hat{q}_i &= \hat{f}_i \text{ on } \Lambda_i,  && \text{(conservation of mass on edge)} \label{eq:2}\\ 
    [[\hat{q}]]_b &= 0 \text{ for } b \in \mathsf{B},  && \text{(conservation of mass at bifurcation)}\label{eq:3}
\end{align}
where $\hat{q}_i$ denotes the cross-section flux along an edge, $\hat{p}_i$ denotes the average pressure on an edge, $\mathcal{R}$ denotes the flow resistance and $\partial_s$ denotes the spatial derivative along the (one-dimensional) edge. Further
\begin{align}
    [[\hat{q}]]_b =  \sum_{i \in \mathsf{E}_{in}(b)} \hat{q}_i(b) - \sum_{i \in \mathsf{E}_{out}(b)} \hat{q}_i(b)
\end{align}
is used to denote the jump in cross-section flux over a bifurcation, with $q= \mathcal{I}_1 q_1+\mathcal{I}_2 q_2+...+\mathcal{I}_n q_n$ denoting the global flux. To close the system we further assume the pressure is continuous over each bifurcation point. 
 
\begin{figure}
\centering
\begin{subfigure}{0.4\textwidth}
    \centering
    \includegraphics[width=0.99\textwidth]{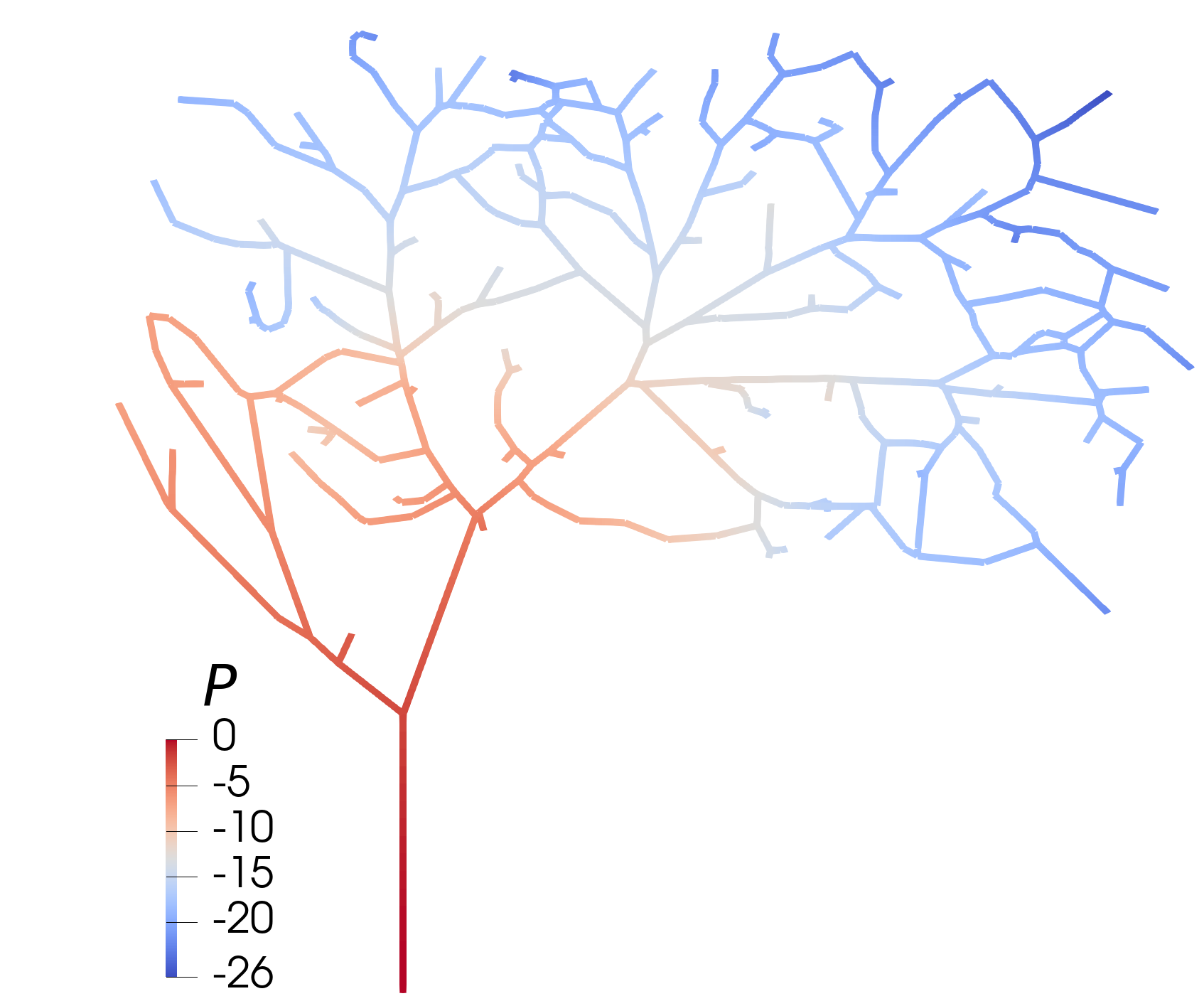}
    \caption{}
    \label{fig:pial-pressure}
\end{subfigure}
\begin{subfigure}{0.4\textwidth}
    \centering
    \includegraphics[width=0.99\textwidth]{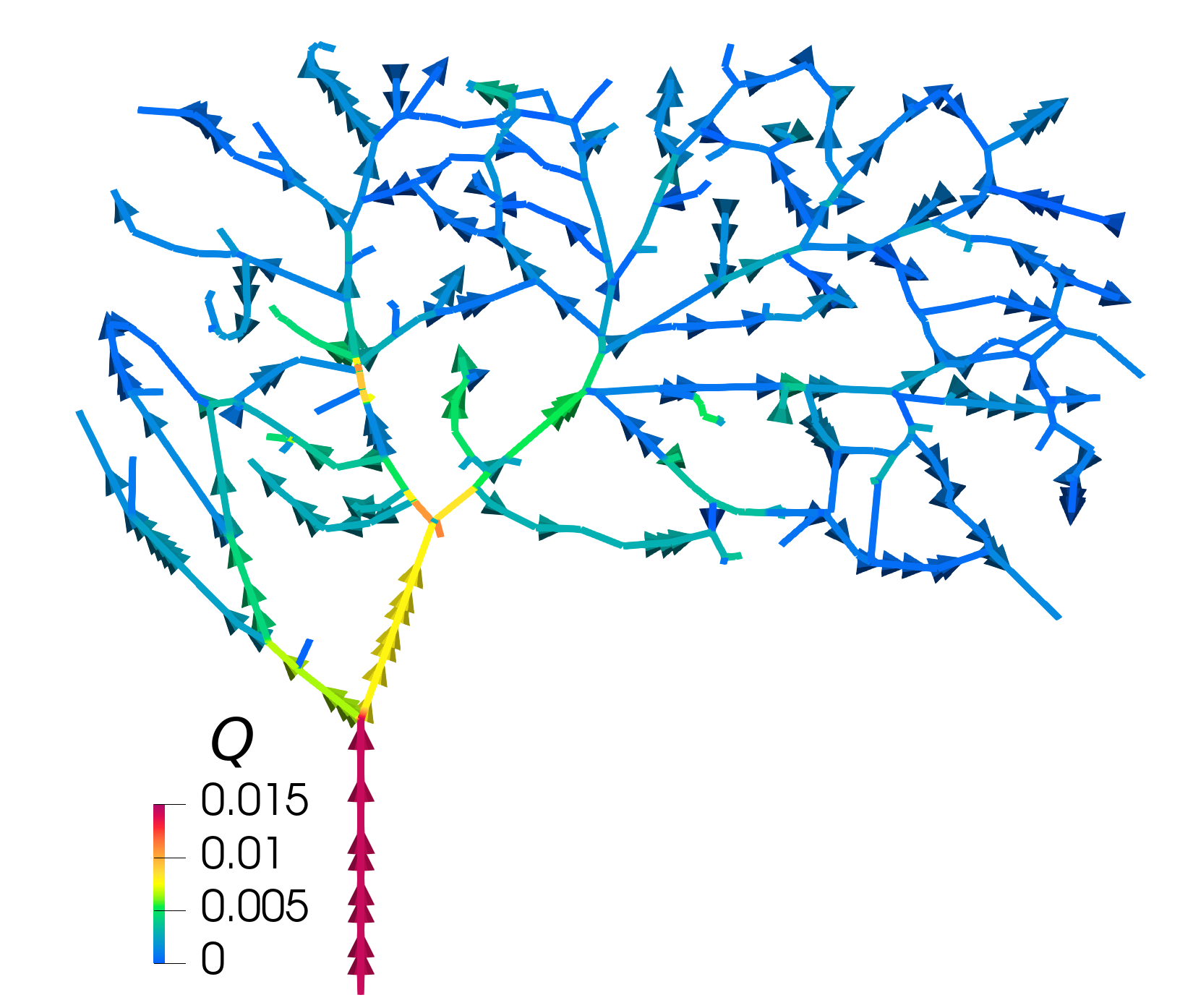}
    \caption{}
    \label{fig:pial-flux}
\end{subfigure}
\label{fig:pial}
\caption{Simulation of fluid flow in the pial blood vessel network of a rodent \cite{topological-kleinfeld}. The network consists of 417 edges and 389 vertices, of which 320 are bifurcation points. The pressure (a) and cross-section flux (b) were computed using the \texttt{HydraulicNetwork} model. A linear pressure drop was ascribed from inlet to outlets.}
\end{figure}

\noindent
Following \cite{notaro2016mixed}, the dual mixed variational form associated with this model reads:  Find $\hat{q} \in V$ and $(\hat{p}, \lambda) \in M$
\begin{align}
    a(\hat{q},\hat{\psi})+b(\psi,(\hat{p}, \lambda)) &= 0 \\
    b((\phi, \xi), \hat{q}) &= (\hat{f}, \phi)
\end{align}
for all $\psi \in V$ and $(\phi, \xi) \in M$, where
\begin{align}
    a(\hat{q}, \psi) &= \sum_{i=1}^n (\mathcal{R} \hat{q}_i, \Psi)_{\Lambda_i}, \\
    b(\psi, (\hat{p}, \lambda)) &= -(\partial_s \psi, \hat{p}) + \sum_{b \in \mathsf{B}} [[\Psi]]_b \lambda_b.
\end{align}
Here $\lambda=(\lambda_1, \lambda_2, \lambda_3, ..., \lambda_m) \in \mathbf{R}^{m}$ is a Lagrange multiplier used to impose conservation of mass at each bifurcation point, and the function spaces are
\begin{align*}
V &= \cup_{i=1}^n H^1(\Lambda_i), \\
M &= L^2(\Lambda) \times \mathbb{R}^m.
\end{align*}

\section*{Software overview}

\subsection*{The \texttt{FenicsGraph} class}

The main component of \textsc{graphnics} is the \texttt{FenicsGraph} class, which inherits from the \textsc{NetworkX} \texttt{DiGraph} class. The \texttt{FenicsGraph} class provides a function for meshing the network; meshfunctions are used to relate the graph structure to the cells and vertices in the mesh. Tangent vectors $\boldsymbol{\tau}_i$ are computed for each edge and stored as edge attributes for the network. This is then used in a convenience class function \textit{dds$\_$i} which returns the spatial derivative $\partial_s f_i = \nabla f_i \cdot \boldsymbol{\tau}_i$ on the edge.

\subsection*{Network models}

\textsc{Graphnics} can be used to create and solve network flow models in \textsc{FEniCS}. A first example of this is shown in the \texttt{NetworkPoisson} model, which can be assemble and solved using standard \textsc{FEniCS} functionality. For discretizations that involve e.g. jump terms, the edge and vertex iterators in \textsc{NetworkX} are used to assemble contributions to the block matrix. This approach is used in the \texttt{HydraulicNetwork} and \texttt{NetworkStokes} flow models. The \texttt{HydraulicNetwork} model implements \eqref{eq:1}-\eqref{eq:3} using the variational formulation from \cite{notaro2016mixed}. The \texttt{NetworkStokes} model solves a Brinkman-Stokes model \cite{daversin2022geometrically, gjerde2023network} which includes an axial viscosity term in the momentum balance equation.

\section*{Demos}

The \textsc{graphnics} library further contains a collection of demos showcasing
\begin{itemize}
\item the essential functionality of \textsc{graphnics}, including how to construct networks, how to mesh the domain and how to solve simple flow models,
\item how \textsc{graphnics} can be used to simulate pulsatile flow in complex biological networks,
\item how \textsc{graphnics} can be used together with \textsc{fenics}$\_$\textsc{ii} to solve coupled 1d-3d flow models.
\end{itemize}

In particular, the pulsatile flow demos extend the simulations from \cite{daversin2022geometrically} to non-trivial networks. There, a network model was presented modelling the flow of Cerebrospinal Fluid (CSF) through so-called Perivascular Spaces (PVSs). The driving forces for this flow has recently attracted extensive research interest due to its possible role in the development of neural diseases.

\subsection*{Using \textsc{graphnics} to simulate PVS flow}
Experimental results indicate that tracers move through the PVS in lockstep with cardiac induced arterial pulsations \cite{mestre2018flow}. This has given rise to the hypothesis that bulk flow of CSF through the brain is driven by arterial wall movement. So far, this has not been replicated in numerical simulations \cite{daversin2020mechanisms, daversin2022geometrically}. 

One caveat of these simulations is that they are performed on a simple network containing only one bifurcation. In  \cite{bedussi2018paravascular}, tracer movement was observed in places where the actual diameter oscillations were negligible. Based on this, the authors suggest that upstream or downstream  arterial wall movements contribute to the overall flow of CSF through the PVS. It is possible that bulk flow of CSF driven by arterial pulsations only occurs in networks of a certain complexity. With this in mind, two of the \textsc{graphnics} demos extend the simulations in \cite{daversin2020mechanisms, daversin2022geometrically} to non-trivial networks.

\subsubsection*{Cardiac pulsations drives purely pulsatile flow in a complex vascular network}

The effect of cardiac wall pulsations was simulated in a vascular network found in a rat carcinoma \cite{secomb1998theoretical}. The PVS was modeled as an annular flow channel with outer radius $R_2=3R_1$, where $R_1$ is the arterial wall radius given in the vasculature dataset. Flow of CSF through the PVS was modelled using the network Stokes model derived in \cite{gjerde2023network}. Wall pulsations were introduced uniformly using experimental data for the relative displacement of the arterial wall from \cite{mestre2018flow}. 

\begin{figure}
    \centering
\begin{subfigure}{0.45\textwidth}
    \includegraphics[width=0.99\textwidth]{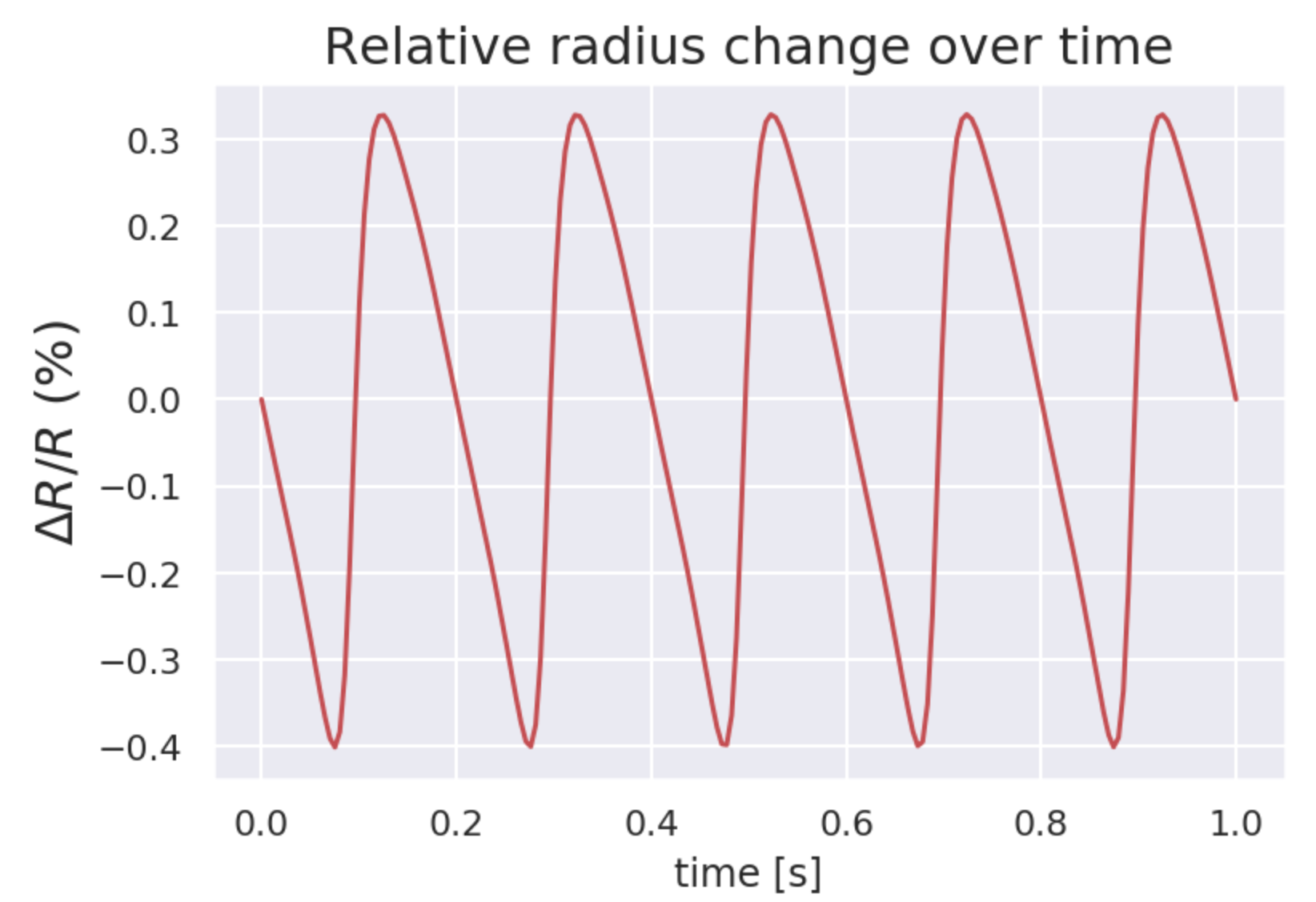}
    \caption{}
    \label{fig:cardiaca}
\end{subfigure}
\begin{subfigure}{0.45\textwidth}
    \includegraphics[width=0.99\textwidth]{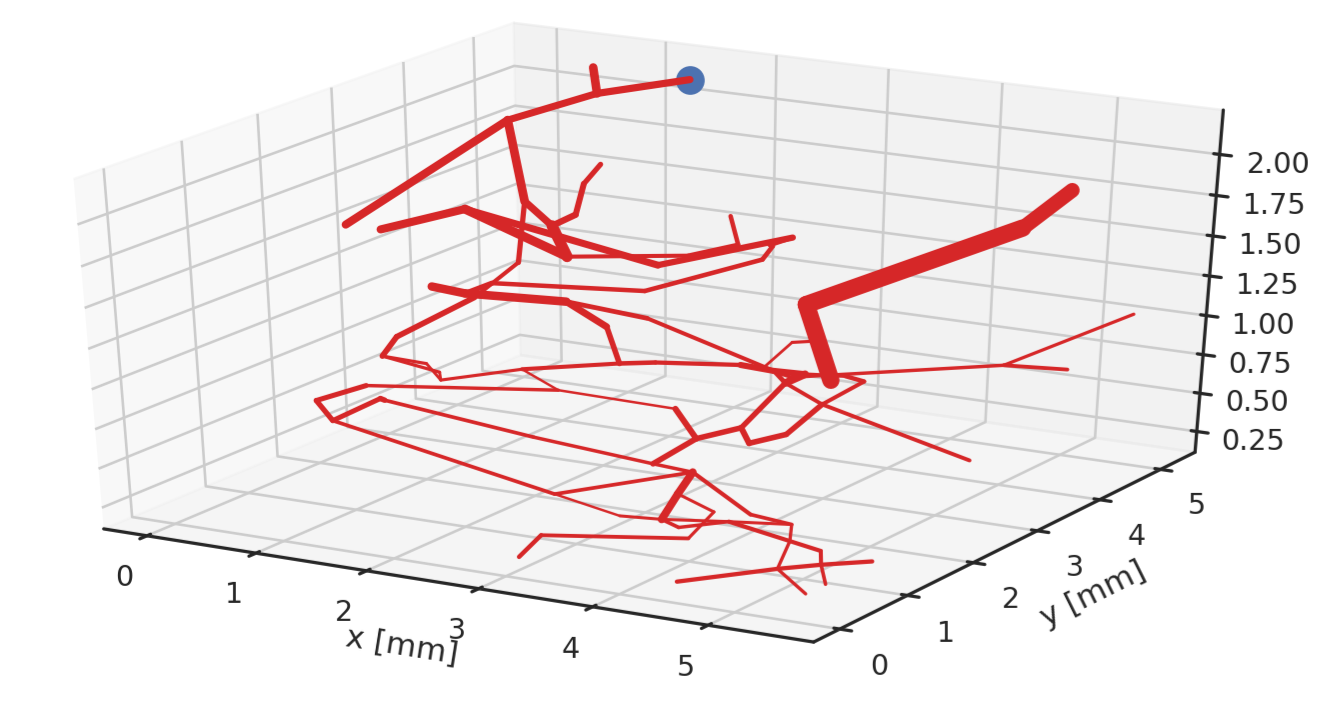}
    \caption{}
        \label{fig:cardiacc}
\end{subfigure}

\vspace{1em}

\begin{subfigure}{0.45\textwidth}
    \includegraphics[width=0.99\textwidth]{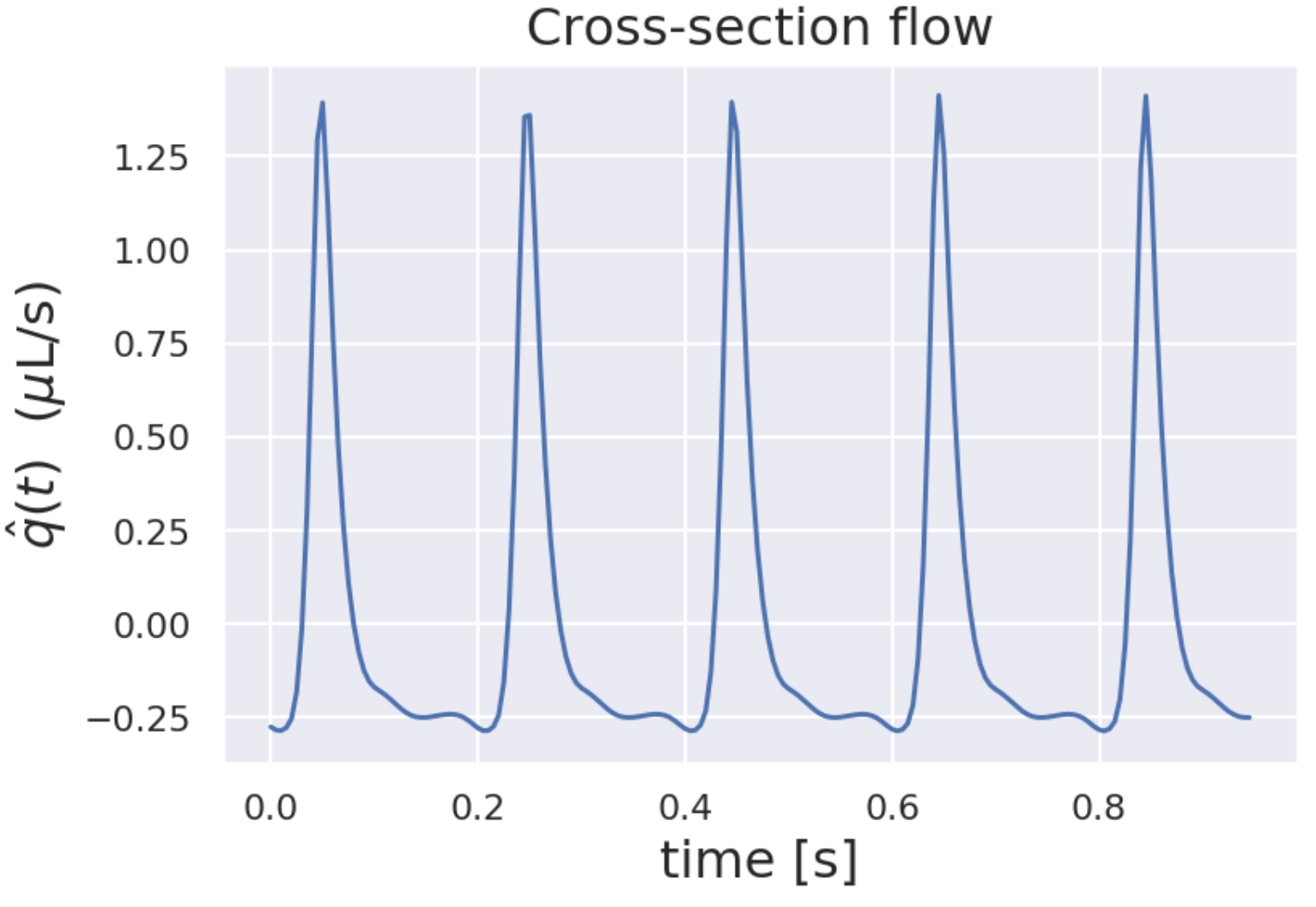}
    \caption{}
        \label{fig:cardiacd}
\end{subfigure}
\begin{subfigure}{0.45\textwidth}
        \includegraphics[width=0.99\textwidth]{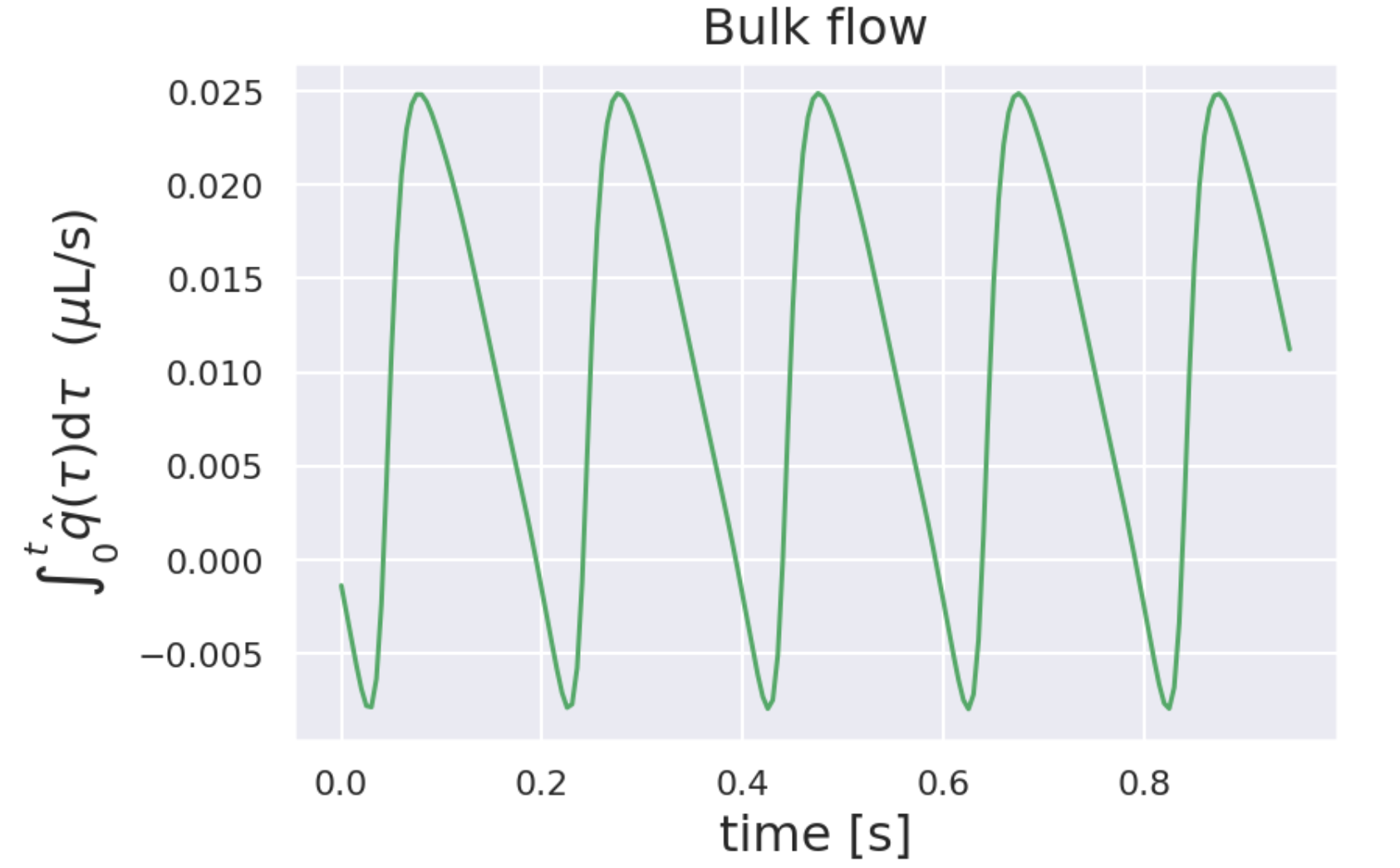}
        \caption{}
            \label{fig:cardiace}
\end{subfigure}

\begin{subfigure}{0.6\textwidth}
    \includegraphics[width=0.99\textwidth]{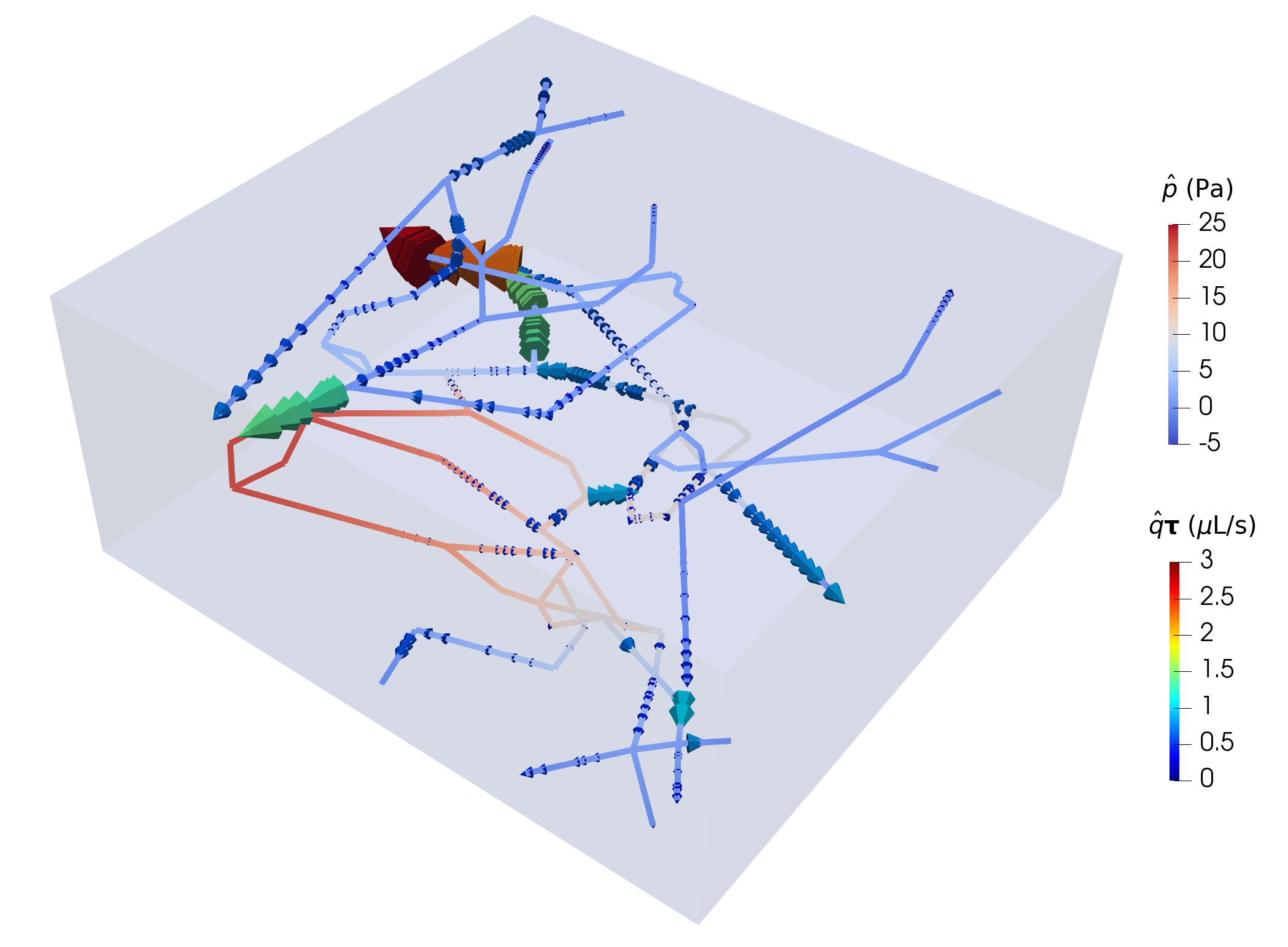}
    \caption{}
        \label{fig:cardiacb}
\end{subfigure}

    \caption{Cardiac wall motion in a vascular network induced oscillatory flow. Pressure and cross-section flux were computed using the  \texttt{NetworkStokes} model. Pulsatile cardiac wall motion was imposed uniformly at all the branches in a vascular network from a rat carcinoma with varying vessel diameter (\ref{fig:cardiacc}). The relative diameter was changed according to (\ref{fig:cardiaca}) The contraction/expansion of the arterial wall drives fluid in an out of the domain, leading to peak outflow rates of 3$\mu$L/s (\ref{fig:cardiacb}). Tracking flow at a specific outlet marked in (\ref{fig:cardiacc}), we see that the resulting flow is pulsatile (\ref{fig:cardiacd}) with no appreciable net flow over time (\ref{fig:cardiace}).}
    \label{fig:cardiac}
\end{figure}

The results are shown in Figure \ref{fig:cardiac}. Figure     \ref{fig:cardiaca} shows the arterial radius as a function of time. A snapshot of the solution at peak outflow is shown in Figure \ref{fig:cardiacb}. As in \cite{daversin2020mechanisms, daversin2022geometrically}, cardiac induced wall movement was found to induce purely oscillatory back-and-forth flow patterns (Figure \ref{fig:cardiacd}), with no appreciable bulk flow (Figure \ref{fig:cardiace}). 

\subsubsection*{Vasomotion drives net flow of CSF in arterial trees of sufficiently complexity}

Recently, it has been proposed that vasomotion might drive bulk flow of CSF through the PVS \cite{van2020vasomotion}. This was previously investigated via simulations in \cite{daversin2022geometrically}; there, vasomotion was modelled as a travelling sinusoidal wave with wave frequency 0.1Hz and wavelength 8mm. No appreciable net flow was found to be induced in an arterial tree with two generations.

\begin{figure}
    \centering
\begin{subfigure}{0.9\textwidth}
        \includegraphics[width=0.99\textwidth]{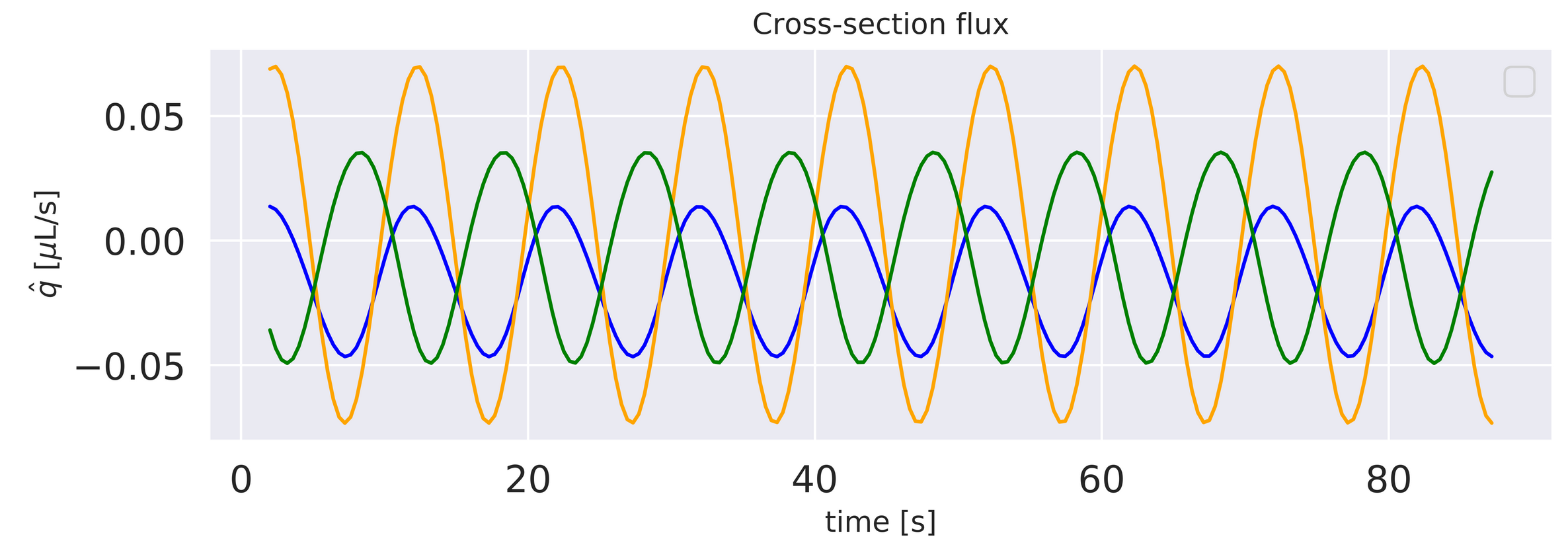}
        \caption{}
        \label{fig:vasomotiona}
\end{subfigure}

\begin{subfigure}{0.25\textwidth}
        \includegraphics[width=0.99\textwidth]{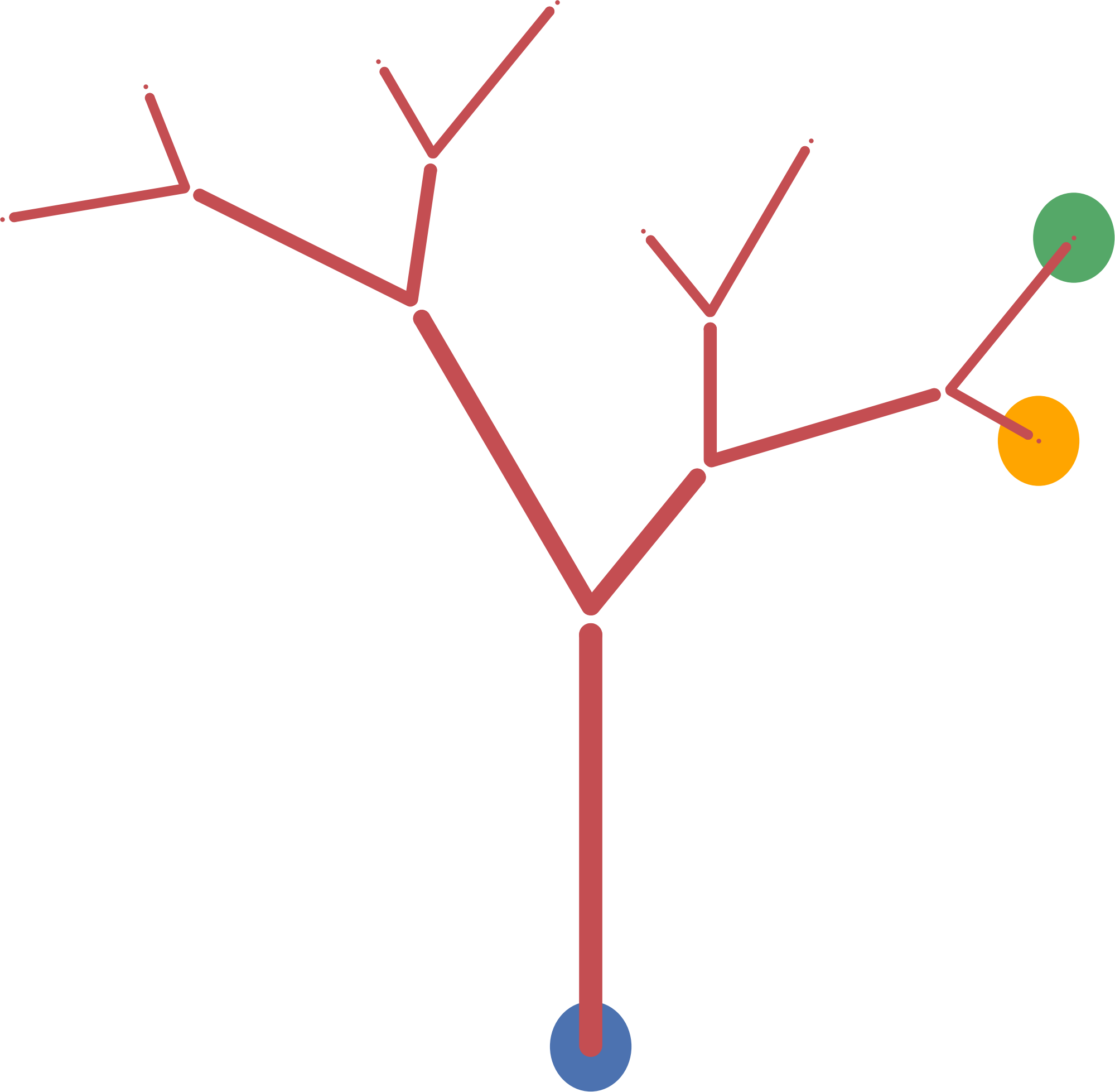}
        \caption{}
        \label{fig:vasomotionb}
\end{subfigure}
\begin{subfigure}{0.5\textwidth}
        \includegraphics[width=0.99\textwidth]{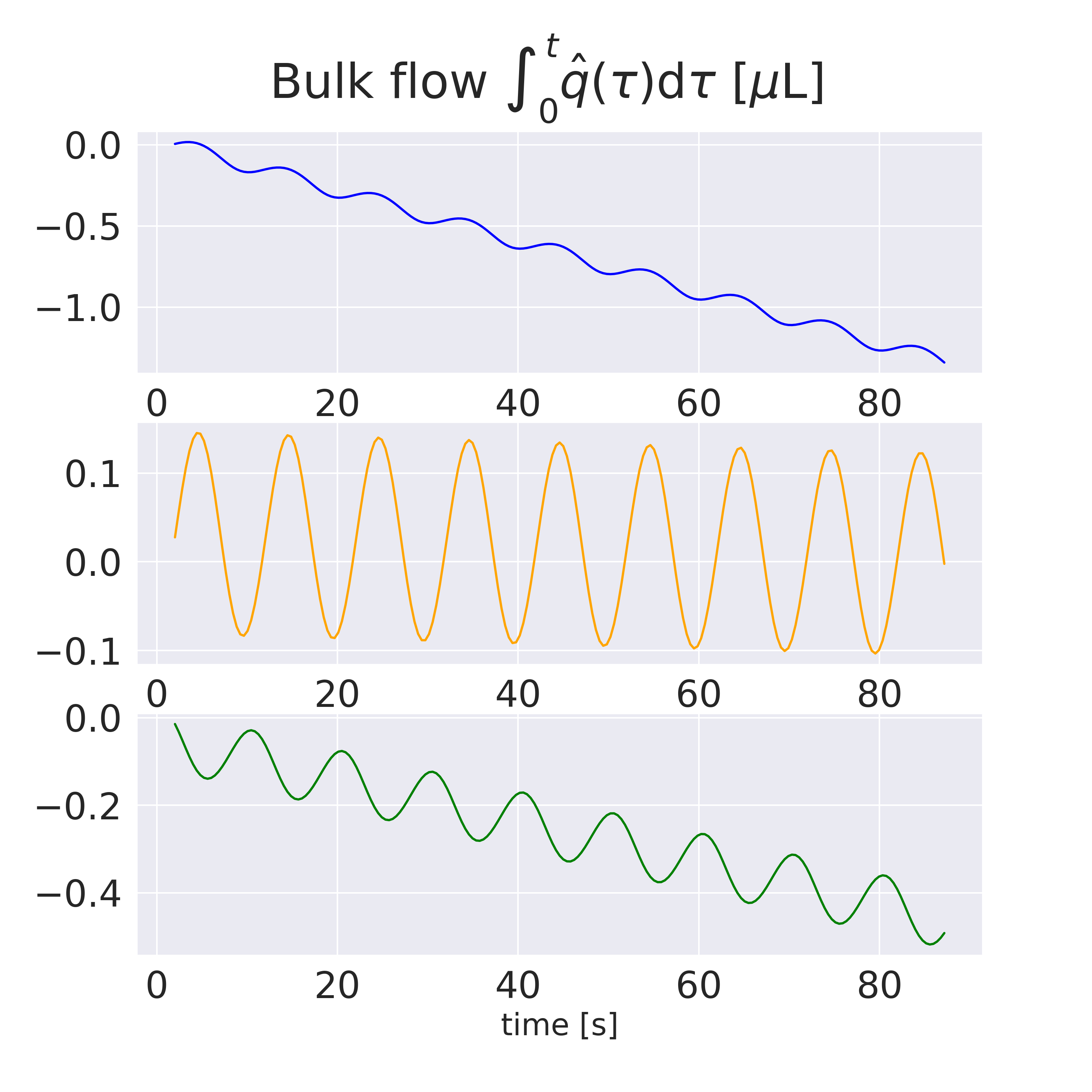}
        \caption{}
        \label{fig:vasomotiond}
\end{subfigure}

\begin{subfigure}{0.99\textwidth}
        \includegraphics[width=0.32\textwidth]{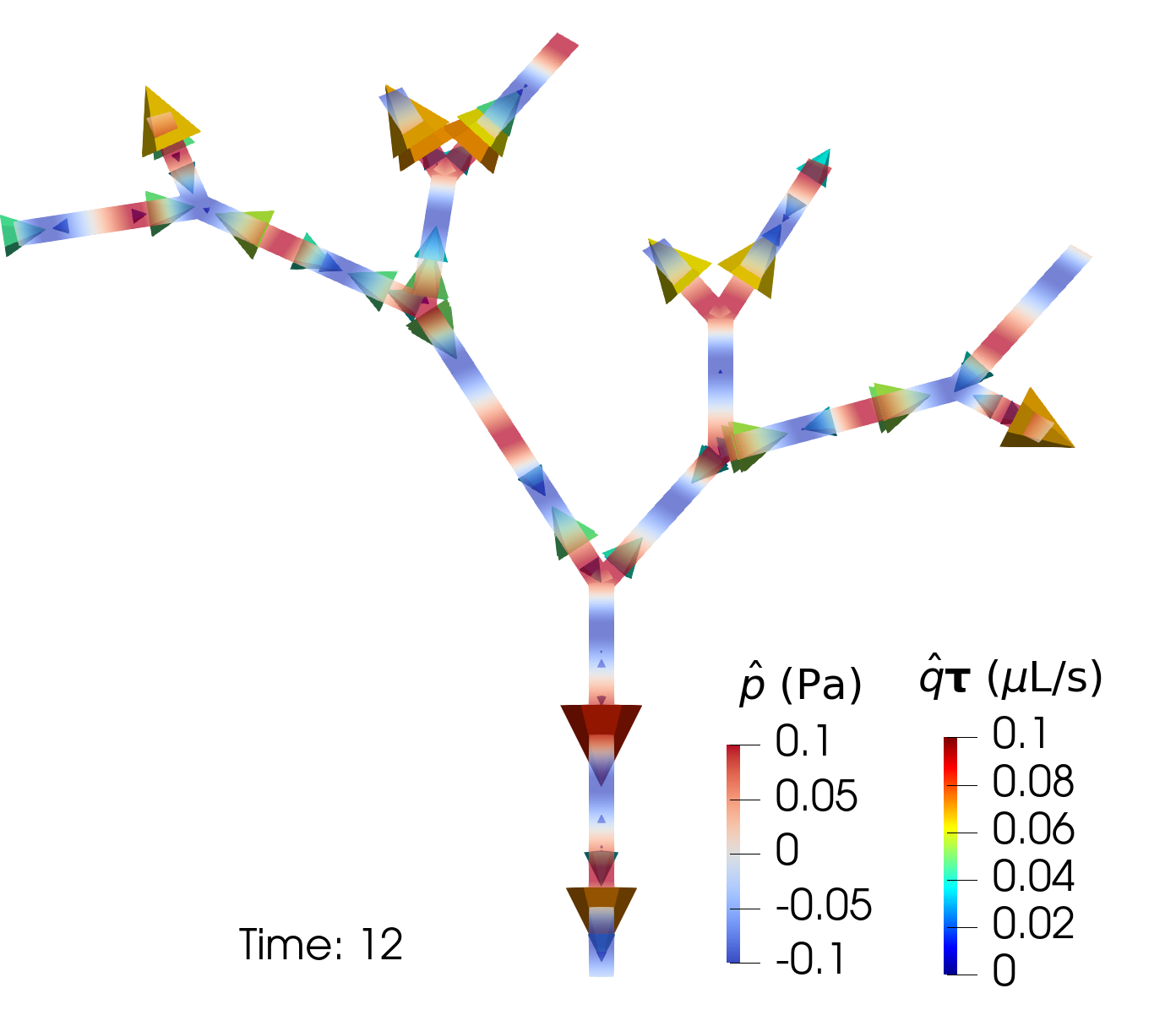}
        \includegraphics[width=0.32\textwidth]{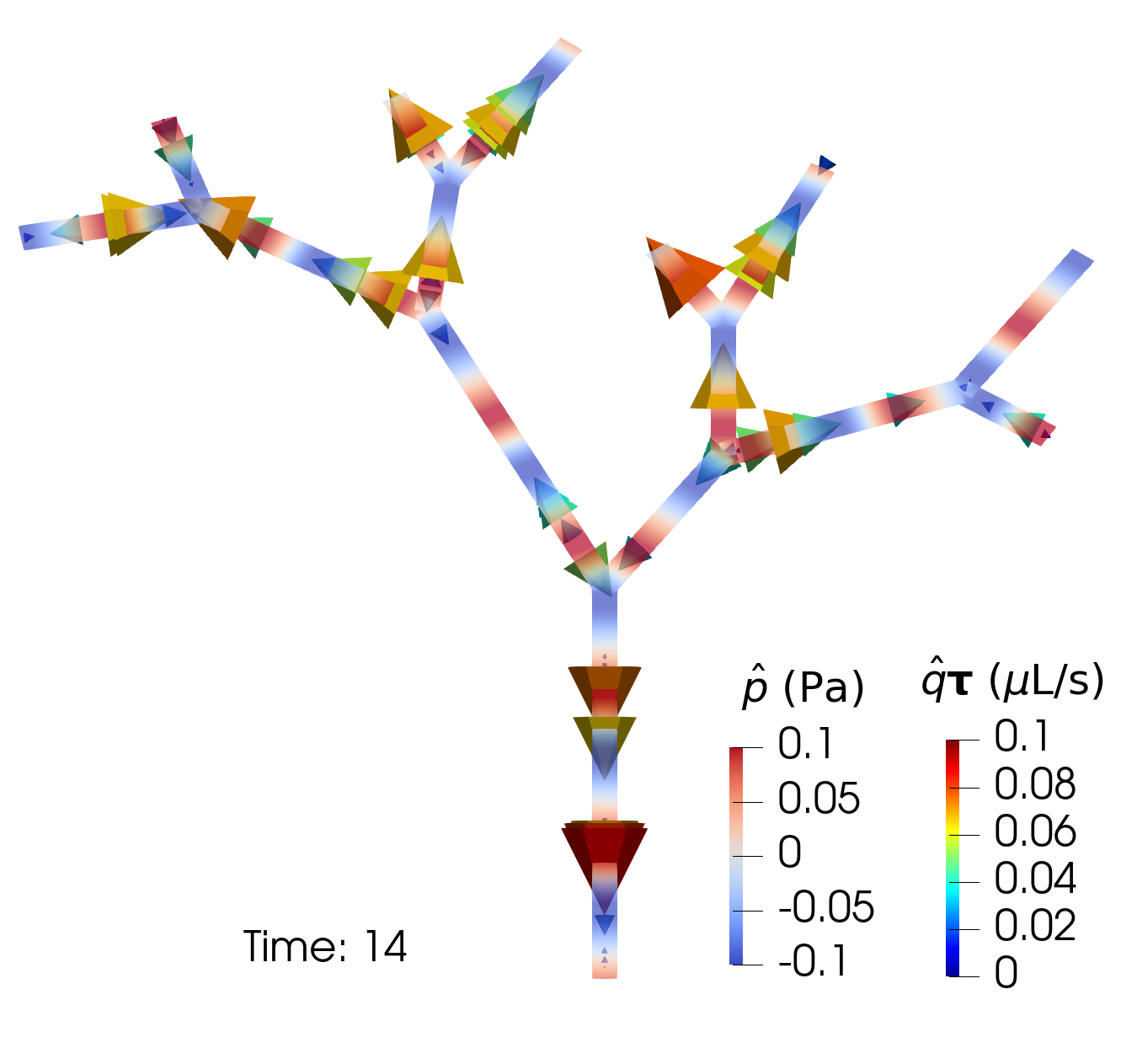}
		\includegraphics[width=0.32\textwidth]{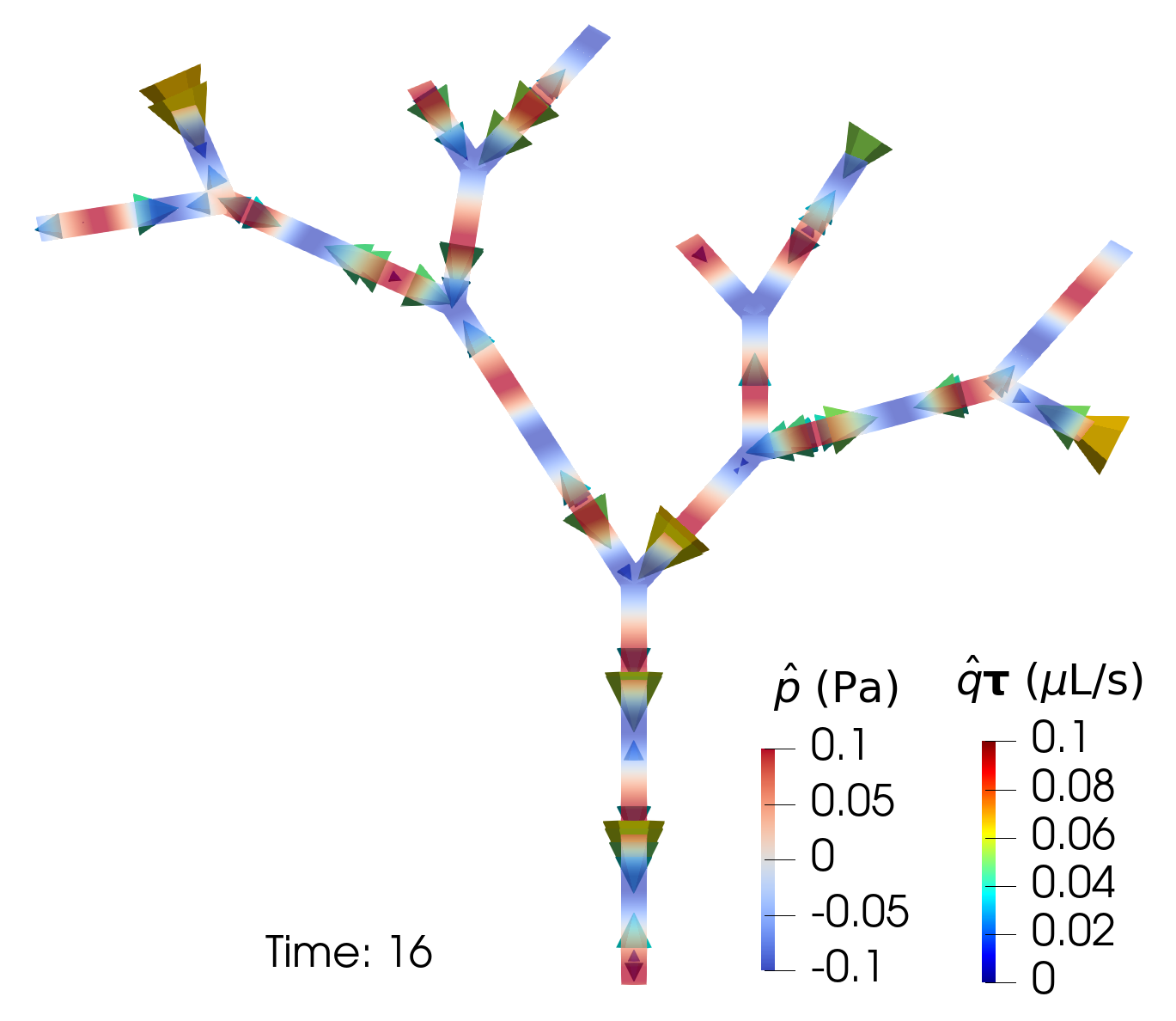}
        \caption{}     
        \label{fig:vasomotionc}
\end{subfigure}

    \caption{Vasomotion induces net flow in an arterial tree with varying radius (\ref{fig:vasomotionb}). Vasomotion was modelled as a travelling sinusoidal wave that expands/contracts the arterial wall 10$\%$ of the initial radius. A snapshot of the solution (\ref{fig:vasomotionc}) shows complex flow patterns. Tracking the flow at the inlet and two outlets shows oscillatory flow (\ref{fig:vasomotiona}), with bulk flow entering the leaf nodes and leaving the root node in the arterial tree (\ref{fig:vasomotiond}).}
    \label{fig:vasomotion}
\end{figure}

Using \textsc{graphnics} to handle more complex networks, it was found that a similar simulation setup produces net flow in an arterial tree with four generations. The arterial tree was generated using \textsc{NetworkGen} \cite{networkgen}, with the radius change at the bifurcations obeying Murray's law. The results are shown in Figure \ref{fig:vasomotion}. A snapshot of the solution is shown in Figure \ref{fig:vasomotionc}. Travelling vasomotion waves were found to introduce complex, oscillatory flow patterns in the arterial tree. Tracking the net flow at the inlets and two selected outlets, it was found that there is a bulk fluid flow moving from leaf nodes to root node (Figure \ref{fig:vasomotiond}).

\section*{Acknowledgments}
We thank Marie Rognes and Barbara Wohlmuth for helpful discussions on network modelling and discretization, Miroslav Kuchta, Cecile Daversin-Catty and Jørgen Dokken for their input on the implementation and Pablo Blinder and David Kleinfeld for sharing data.

\bibliographystyle{plain}


\end{document}